\documentclass[final, 3p, 12pt]{elsarticle}
\usepackage{amsthm}
\usepackage{amsmath}
\usepackage{amssymb}
\usepackage{setspace}
\usepackage{simplemargins}
 \setallmargins{1.0in}
\usepackage{rotating}
\usepackage{epsfig}
\usepackage{float}
\usepackage{relsize}
\usepackage{verbatim}
\usepackage{amsfonts}
\usepackage{fancyhdr}

\newcommand\+{\mkern2mu}
\usepackage{longtable}
\usepackage{endnotes}
\let\footnote=\endnote

  \usepackage{scalerel,stackengine}
\stackMath
\newcommand\reallywidehat[1]{%
\savestack{\tmpbox}{\stretchto{%
  \scaleto{%
    \scalerel*[\widthof{\ensuremath{#1}}]{\kern-.6pt\bigwedge\kern-.6pt}%
    {\rule[-\textheight/2]{1ex}{\textheight}}
  }{\textheight}%
}{0.5ex}}%
\stackon[1pt]{#1}{\tmpbox}%
}
\makeatletter
\newcommand*\rel@kern[1]{\kern#1\dimexpr\macc@kerna}
\newcommand*\widebar[1]{%
  \begingroup
  \def\mathaccent##1##2{%
    \rel@kern{0.8}%
    \overline{\rel@kern{-0.8}\macc@nucleus\rel@kern{0.2}}%
    \rel@kern{-0.2}%
  }%
  \macc@depth\@ne
  \let\math@bgroup\@empty \let\math@egroup\macc@set@skewchar
  \mathsurround\z@ \frozen@everymath{\mathgroup\macc@group\relax}%
  \macc@set@skewchar\relax
  \let\mathaccentV\macc@nested@a
  \macc@nested@a\relax111{#1}%
  \endgroup
}
\makeatother
\theoremstyle{definition}
\newtheorem*{pf}{Proof}
\theoremstyle{theorem}
\newcounter{dummyt} \numberwithin{dummyt}{section}
\newcounter{dummym} \numberwithin{dummym}{section}
\newcounter{dummya} \numberwithin{dummya}{section}
\newcounter{dummyc} \numberwithin{dummyc}{section}
\newcounter{dummyl} \numberwithin{dummyl}{section}
\newcounter{dummyd} \numberwithin{dummyd}{section}
\newcounter{dummye} \numberwithin{dummye}{section}
\newtheorem{thm}[dummyt]{Theorem}
\newtheorem{meth}[dummym]{Method}

\newtheorem{lem}[dummyl]{Lemma}
\theoremstyle{definition}
\newtheorem{defn}[dummyd]{Definition}
\newtheorem{ex}[dummye]{Example}
\usepackage{url}
\usepackage{xcolor}
\begin{document}
\newcommand{\oankm}{\text{OA}(N,k-1,s,t)}
\newcommand{\oan}{\text{OA}(N,k,s,t)}\newcommand{\oano}{\text{OA}(N,k_0,s,t)}
\newcommand{\oant}{\text{OA}(N,k,2,t)}
\newcommand*{\tabref}[1]{\tablename~\ref{#1}}
\newcommand{\am}{\boldsymbol{A}}
\newcommand{\av}{\boldsymbol{a}}
\newcommand{\bv}{\boldsymbol{b}}
\newcommand{\cv}{\boldsymbol{c}}
\newcommand{\dv}{\boldsymbol{d}}

\newcommand{\fnst}{f(N,s,t)}
\newcommand{\fnt}{f(N,2,t)}
\newcommand{\ftnt}{f(2N,2,t+1)}
\newcommand{\fcan}{f^{\text{CA}}_{\lambda}(N,s,t)}
\newcommand{\fpan}{f^{\text{PA}}_{\lambda}(N,s,t)}
\newcommand{\fcant}{f^{\text{CA}}_{\lambda}(N,2,t)}
\newcommand{\fpant}{f^{\text{PA}}_{\lambda}(N,2,t)}

\newcommand{\kin}{k_{\mathrm{input}}}
\newcommand{\kout}{k_{\mathrm{out}}}

\newcommand{\fxoptt}{f(\X^{opt}_{N})}
\newcommand{\yoptt}{\Y^{opt}_{N,k}}
\newcommand{\xoptt}{\X^{opt}_{N,k}}
\newcommand{\x}{\boldsymbol{x}}
\newcommand{\z}{\boldsymbol{z}}
\newcommand{\M}{\boldsymbol{M}}
\newcommand{\N}{\boldsymbol{N}}
\newcommand{\ZZ}{\boldsymbol{Z}}
\newcommand{\HH}{\boldsymbol{H}}
\newcommand{\I}{\boldsymbol{I}}
\newcommand{\1}{\boldsymbol{1}}
\newcommand{\zz}{\boldsymbol{0}}
\newcommand{\aaa}{\boldsymbol{a}}
\newcommand{\bb}{\boldsymbol{b}}
\newcommand{\B}{\boldsymbol{B}}
\newcommand{\cc}{\boldsymbol{c}}
\newcommand{\C}{\boldsymbol{C}}
\newcommand{\CC}{\mathcal{C}}
\newcommand{\uu}{\boldsymbol{u}}
\newcommand{\vv}{\boldsymbol{v}}
\newcommand{\rrr}{\boldsymbol{r}}
\newcommand{\xx}{\boldsymbol{x}}
\newcommand{\y}{\boldsymbol{y}}
\newcommand{\dd}{\boldsymbol{d}}
\newcommand{\vvec}{\boldsymbol{vec}}
\newcommand{\Y}{\boldsymbol{Y}}
\newcommand{\X}{\boldsymbol{X}}
\newcommand{\D}{\boldsymbol{D}}
\newcommand{\A}{\boldsymbol{A}}
\newcommand{\tA}{\tilde{\boldsymbol{A}}}
\newcommand{\XX}{\boldsymbol{X}^{\top}\boldsymbol{X}}
\newcommand{\XT}{\boldsymbol{X}^{\top}}
\newcommand{\bs}{\boldsymbol}
\newcommand{\orr}{\text{\ or \ }}
\newcommand{\md}[2][4]{#2\ \ (\mathrm{mod}\ \ #1)}
\newcommand{\rr}{\mbox{I\!\!\,R}}
\newcommand{\es}{\mathrm{E}(s^2)}
\newcommand{\smax}{s_{\mathrm{max}}}
\newcommand{\freq}{f_{\smax}}
\newcommand{\dfk}[1]{\Delta f^k(#1)}
\newcommand{\Z}{\mathbb{Z}}
\newcommand{\mth}[1]{\multicolumn{3}{c|}{#1}}
\newcommand{\mtw}[1]{\multicolumn{2}{c|}{#1}}
\newcommand{\GWP}{\mathrm{GWP}}
\begin{frontmatter}
\title{Integer programming for classifying orthogonal arrays}
\author[AFIT]{Dursun A.~Bulutoglu\corref{cor1}}
\ead{dursun.bulutoglu@gmail.com}
\author[WVU]{Kenneth J.~Ryan}
\ead{kjryan@mail.wvu.edu}
\address[AFIT]{Department of Mathematics and Statistics, Air Force Institute of Technology,\\Wright-Patterson Air Force Base, Ohio 45433, USA}
\address[WVU]{Department of Statistics, West Virginia University,\\ Morgantown, West Virginia 26506, USA}
\cortext[cor1]{Corresponding author}
\journal{Australasian Journal of Combinatorics}
\begin{abstract}
Classifying orthogonal arrays (OAs) is a well-known important class of problems that asks for finding all non-isomorphic, non-negative integer solutions
 to a class of systems of constraints. Solved instances are scarce.
We develop two new methods based on finding all non-isomorphic solutions of two  novel integer linear programming formulations for classifying all non-isomorphic $\oan$ given a set of all non-isomorphic $\oankm$. 
 We also establish the concept of orthogonal design  equivalence
 (OD-equivalence) of $\oant$ to reduce the number of integer linear programs (ILPs) whose all non-isomorphic  
  solutions need to be enumerated by our methods. 
 For each ILP, we determine the largest group of permutations that can be exploited with the branch-and-bound (B\&B) with isomorphism pruning 
algorithm of Margot [Discrete Optim.~4 (2007), 40-62]  without losing isomorphism classes of $\oant$.
Our contributions bring the classifications of all
non-isomorphic $\text{OA}(160,k,2,4)$ for $k=9,10$ and $\text{OA}(176,k,2,4)$ for $k=5,6,7,8,9,10$ within computational reach.
These are the smallest $s=2$, $t=4$ cases for which classification results are not available in the literature.
\end{abstract}
\end{frontmatter}
\section{Introduction}
This  work develops methods based on state-of-the-art discrete optimization algorithms to solve several instances of
 a fundamental class of problems in combinatorics. We establish  novel  methods to classify orthogonal arrays (OAs); see Definition~\ref{def:OA}.  These methods are based on the  
 branch-and-bound (B\&B) with isomorphism pruning algorithm~\cite{Margot07}  or  a basic B\&B algorithm implemented in {\tt CPLEX12.5.1} interactive  and  our newly formulated orthogonal design equivalence (OD-equivalence) of $2$-symbol OAs. 

\begin{defn}\label{def:OA}
 An $N\times k$ array $\Y$ whose entries are symbols from $\{0,\ldots,s-1\}$ is an OA of strength $t$, $1
\leq t \leq k$, denoted by $\oan$, if each of the $s^t$ symbol
combinations from $\{0,\ldots,s-1\}^t$ appears $\lambda=N/s^t$ times in every $N \times  t $ subarray of  $\Y$.
 \end{defn}
Two $N \times k$ arrays with entries  from $\{0,\ldots,s-1\}$   are
called {\em isomorphic} if one can be obtained from the other by  applying a sequence of permutations to the  columns, the rows, and to  the elements $\{0,\ldots,s-1\}$ within each column~\cite{Stufken}. 
Any  $N\times k$ array that is isomorphic to an $\oan$ is also an $\oan$. Hence it suffices to classify all non-isomorphic $\oan$ to find all $\oan$. Classifying all non-isomorphic $\oan$ is an established problem~\cite{BulutogluM,Schoen,Stufken}, and classification results are scarce~\cite{Stufken}. 
For applications of OAs in statistics and their connection to error correcting codes, see~\cite{Hedayat}. 

Throughout this paper,  an $s$-symbol design refers to an $N\times k$ array whose entries are symbols from $\{0,\ldots,s-1\}$ 
except in Sections~\ref{sec:fold} and~\ref{sec:validation}, where
 $2$-symbol OAs or  designs are recoded with entries from $\{-1, 1\}$ for OD-equivalence reductions.
 

Bulutoglu and Margot~\cite{BulutogluM} found all non-isomorphic solutions to one integer linear program (ILP) to directly find a set of all non-isomorphic $\oan$ for many $N,k,s,t$ combinations.
Schoen, Eendebak, and Nguyen~\cite{Schoen} developed the Minimum Complete Set (MCS) algorithm for enumerating OAs up to isomorphism. MCS is a backtrack algorithm with isomorph rejection. Schoen, Eendebak, and Nguyen~\cite{Schoen} also used their MCS algorithm to classify all non-isomorphic $\oan$ for many $N,k,s,t$ combinations--including all those in Bulutoglu and Margot~\cite{BulutogluM}--and mixed level OAs where all columns do not have the same number of distinct symbols. For $\lambda>1$, the MCS algorithm is up to $3$ orders of magnitude faster than the other algorithms in the literature.
For $\lambda=1$,  Kokkala and \"{O}sterg\r{a}rd~\cite{Kokkala15} recently obtained a set of all non-isomorphic $\text{OA}(q^2,4,q,2)$ from a set of all non-isomorphic $\text{OA}(q^2,3,q,2)$ for $3 \leq q\leq 8$ by solving exact cover problems, and Egan and Wanless~\cite{Egan} enumerated all non-isomorphic  $\text{OA}(q^2,k,q,2)$ for $q\leq 9$ and $k \leq q+1$ with a faster and more specialized algorithm.
However, they failed to determine whether an OA$(100,5,10,2)$ exists.  
Appa, Magos, and Mourtos~\cite{AMM} also tried, but failed to solve this problem using an algorithm that switches between constraint programming and integer programming search trees, where the inherent symmetries were removed by introducing constraints. Kokkala and \"{O}sterg\r{a}rd~\cite{Kokkala15}  developed an algorithm based on $k$-seeds for classifying all non-isomorphic $\text{OA}(q^3,5,q,3)$ for $q=5,7,8$. Their algorithm first finds all non-isomorphic $k$-seeds and then it completes the resulting $k$-seeds to $\text{OA}(q^3,5,q,3)$. The completion problem is cast as a maximum clique finding problem in a certain graph. Removal of isomorphic designs in  the classifications above were done with \texttt{nauty}~\cite{McKayP} after mapping designs to graphs. A pair of these graphs are isomorphic if and only if their corresponding designs are isomorphic. Our methods in this paper are not competitive with those in~\cite{Egan,Kokkala15} for these $\lambda=1$ cases.

Appa, Magos, and Mourtos~\cite{AMMOLS} and Appa, Mourtos, and Magos~\cite{AMM}  compared integer programming and constraint programming in determining  whether  an 
OA$(q^2,4,q,2)$ exists for $3 \leq q \leq 12$ after removing symmetries by introducing constraints.
They observed that  solving LP relaxations on the nodes of a search tree becomes more beneficial at improving solution times as the number of variables grows.
Our Section~\ref{sec:results} results corroborate their observation.

Stufken and Tang~\cite{Stufken} classified all non-isomorphic OA$(\lambda2^t,t+2,2,t)$, and  Fujii, Namikawa and Yamamoto~\cite{Fuji} classified all non-isomorphic OA$(2^{t+2},t+3,2,t)$ for all $\lambda>0$, $t>0$ with $\lambda,t \in \mathbb{Z}$.
 These are the only infinite classes of OAs for which classification results are available. Hence, it is natural to classify 
more open cases of $2$-symbol OAs. 

In Section~\ref{sec:fold}, we introduce the concept of  our newly formulated OD-equivalence for $2$-symbol designs.
In Section~\ref{sec:isop}, the notions of symmetric ILP and   isomorphic solutions of a symmetric  ILP are introduced. Also, the B\&B  with isomorphism pruning algorithm~\cite{Margot07} for finding a set of all non-isomorphic solutions to an ILP is discussed.

Two new  ILP formulations are introduced in Section~\ref{sec:formulations}. Both of these ILPs cast the constraint satisfaction problem resulting from  adding a new column to an OA$(N,k-1,s,t)$ to get an OA$(N,k,s,t)$. The newly introduced methods  of  Section~\ref{sec:enum},  ILP-extension and Hybrid, use these ILP formulations to obtain a set of all non-isomorphic $\oan$ from a set of all non-isomorphic OA$(N,k-1,s,t)$ by adding columns.
The newly established classification results  of all non-isomorphic $\text{OA}(160,k,2,4)$ for $k=9,10$ and $\text{OA}(176,k,2,4)$ for $k=5,6,7,8,9,10$ are in  
Table~\ref{tab:enum}. These were the smallest unknown $2$-symbol strength $4$ cases. Our results were made computationally feasible by combining: (i) two novel methods in Section~\ref{sec:enum} based on a basic B\&B or B\&B with isomorphism pruning~\cite{Margot07} to find all non-isomorphic or all solutions of our two newly formulated ILPs,  (ii) the concept of OD-equivalence of $2$-symbol  designs in Section~\ref{sec:fold}, and (iii)  exploiting the largest group with
 B\&B with isomorphism pruning~\cite{Margot07}
 that does not lose isomorphism classes of OAs in each ILP.


The theorem within Section~\ref{sec:fold}  made some of the classifications presented in Section~\ref{sec:results} computationally feasible. It justified the use of OD-equivalence of $2$-symbol  OAs to decrease the number of ILPs needed to be solved to find a set of all non-isomorphic $\oant$ by adding columns to a set of  non-isomorphic OA$(N,k-1,2,t)$. 
In Section~\ref{sec:results}, we compare our methods to the MCS 
algorithm~\cite{Schoen} and our ILP-based implementation of McKay's technique~\cite{McKayfree,Radiz} for exploiting symmetries in obtaining larger combinatorial objects from smaller ones to classify $2$-symbol OAs. 
The MCS algorithm cannot classify OA$(160,9,2,4)$ and OA$(176,9,2,4)$ even after it is combined with our newly introduced OD-equivalence reduction of $2$-symbol  OAs. On the other hand, we can classify OA$(160,9,2,4)$ and OA$(176,9,2,4)$ and prove that no OA$(N,10,2,4)$ exists for $N=160,176$  by using a basic B\&B algorithm implemented in {\tt CPLEX12.5.1} combined with McKay's technique and OD-equivalence reductions. However, for classifying OA$(160,8,2,4)$ and OA$(176,8,2,4)$ in which symmetry exploitation is necessary, our methods are $11.44$ and $1.44$ times faster. For cases in which symmetry exploitation is not necessary McKay's technique does not decrease solution times.

The orbit-stabilizer theorem is used in Section~\ref{sec:validation} to derive two formulas that independently count all $\oant$ up to permutations of rows. Since these counts match for each newly classified $\text{OA}(N,k,2,4)$ case, the classification results of Section~\ref{sec:results} are validated. Section~\ref{sec:sum} contains some follow-up discussion including challenge problems and  future research directions. 
\section{OD-equivalence for classifying {\boldmath $\oant$}}\label{sec:fold}
 \setcounter{dummyt}{0} 
\setcounter{dummym}{0} 
\setcounter{dummya}{0} 
\setcounter{dummyc}{0}
\setcounter{dummyl}{0}
\setcounter{dummyd}{0}
\setcounter{dummye}{0}
Throughout this section, $2$-symbol designs or OAs are assumed to have their entries from $\{-1, 1\}$.
First, we define the concept of Hadamard equivalence in~\cite{McKay} and our newly formulated concept of OD-equivalence of $2$-symbol  designs.
 These concepts are used to drastically speed up the Basic Extension method in Section~\ref{sec:enum}  when $s=2$.
\begin{defn}\label{Hadamardeq}
Designs $\Y_1$ and $\Y_2$ are called \emph{Hadamard equivalent} if $\Y_2$ can be obtained from $\Y_1$ by applying signed permutations (permutations that may or may not be followed by sign changes) to the columns or rows of $\Y_1$.
\end{defn}
\begin{defn}\label{ODeq}
 Designs $\X_1$ and $\X_2$ are called \emph{OD-equivalent} if $\left[\1 , \X_1\right]$ and $\left[\1 , \X_2\right]$ are Hadamard equivalent, where $\1$ is the all $1$s column.
\end{defn}

Isomorphic $\oant$ are obviously OD-equivalent. However, there are non-isomorphic $\oant$ that are OD-equivalent. To see this, consider the following example.
\begin{ex}\label{ex:ex}
Let
\begin{equation*}{
[\1 , \X_1]=\left[
\begin{array}{rrrrr}
1 &  -1 & -1 & -1 & -1 \\
1 &  -1 & -1 &  1 &  1 \\
1 &  -1 &  1 & -1 &  1 \\
1 &  -1 &  1 &  1 & -1 \\
1 &   1 & -1 & -1 &  1 \\
1 &   1 & -1 &  1 & -1 \\
1 &   1 &  1 & -1 & -1 \\
1 &   1 &  1 &  1 &  1
\end{array}
\right]\!,\ \
[\1 , \X_2]=\left[
\begin{array}{rrrrr}
1 &  -1 & -1 & -1 &  1 \\
1 &  -1 & -1 &  1 & -1 \\
1 &  -1 &  1 & -1 & -1 \\
1 &  -1 &  1 &  1 &  1 \\
1 &   1 & -1 & -1 &  1 \\
1 &   1 & -1 &  1 & -1 \\
1 &   1 &  1 & -1 & -1 \\
1 &   1 &  1 &  1 &  1
\end{array}
\right]\!,}
\end{equation*}
and  $\av \odot \A$ denote the matrix with the same dimension as $\A$ whose $i$th column is the Hadamard product of the $i$th column of $\A$ and column vector $\av$, where the Hadamard product of vectors $\uu,\vv\in \mathbb{R}^n$ is $\uu\odot\vv=(u_1v_1,u_2v_2,\ldots,u_nv_n)^{\top}$. Designs $\X_i=[\av,\bv,\cv,\dv_i]$ have the following basic structure. The first $3$ columns $[\av,\bv,\cv]$ constitute all $2^3$ combinations $(\pm 1,\pm 1,\pm 1)$, and the last column $\dv_i$ is generated from the defining relations: $\1=\av\odot \bv \odot \cv \odot \dv_1$ and $\1=\bv \odot \cv \odot \dv_2$.
 It is easy to see that $\X_i$ is an $\text{OA}(8,4,2,t_i)$ with $t_i=4-i$ but not $t_i=5-i$ for $i=1,2$.  Although $\X_1$ and $\X_2$ are clearly non-isomorphic, they are OD-equivalent because  $\dv_1=\av \odot \dv_2$ implies $[\1 , \X_1]=\av\odot[\1 , \X_2]$ up to column and row permutations.
\end{ex}
  The generalization of Example~\ref{ex:ex} in Theorem~\ref{thmOD} provides a data compression and makes more classifications possible. The following two definitions and lemma are needed to formulate and prove Theorem~\ref{thmOD}. 

\begin{defn}
Let $\ZZ_1$ and $\ZZ_2$ be OA$(N,k_1,2,t)$ and OA$(N,k_2,2,t)$ with $k_1<k_2$. If $\ZZ_2$ can be obtained by adding $k_2-k_1$ columns to $\ZZ_1$, then  $\ZZ_1$ extends to $\ZZ_2$, or $\ZZ_2$ extends $\ZZ_1$, or 
$\ZZ_2$ is an extension of $\ZZ_1$.
\end{defn}

\begin{defn}[Deng and Tang~\cite{Deng}] 
Let $\Y=[y_{ij}]$ be a $2$-symbol design with   entries of
$\pm 1$ having $N$ rows and $k$ columns, and let
$\ell=\{i_1,i_2,\ldots,i_r\}\subseteq \{1,\ldots,k\}$ be a
non-empty subset of $r$ columns.  Then the values 
\begin{equation*}
J_r(\ell)(\Y):=\sum_{i=1}^N{\prod_{j \in \ell}{y_{ij}}}
\end{equation*}
are called the $J$-characteristics of $\Y$.
\end{defn}
\begin{lem}[Stufken and Tang~\cite{Stufken}]
A $2$-symbol design $\Y$ is an $\oant$ if and only if  $J_r(\ell)(\Y)=0$ for all $\ell \subseteq \{1,\ldots,k\}$ such that $|\ell|=r$ and $r\in \{1,\ldots,t\}$.
 \end{lem}
\begin{thm}\label{thmOD}
Let $\Y$ be an $\text{OA}(N,k_1,2,t)$ that is not an $\text{OA}(N,k_1,2,t+1)$. Let $\Y$ extend to an $\text{OA}(N,k_2,2,t)$. Also, let $\ \ZZ$ be OD-equivalent, but not isomorphic to $\Y$.
 \begin{enumerate}
 \item[(i)] If $t$ is even, then $\ZZ$ is an  $\text{OA}(N,k_1,2,t)$ that extends to an $\text{OA}(N,k_2,2,t)$.
\item[(ii)] If $t$ is odd, then $\ZZ$ is an  OA($N,k_1,2,t-1$), but not necessarily an $\text{OA}(N,k_1,2,t)$. Moreover, $\ZZ$ extends to an $\text{OA}(N,k_2,2,t-1)$.
\end{enumerate}
\end{thm}
\begin{pf}
Let $\Y=[\y_1, \y_2, \cdots, \y_{k_1}]$ and $\ZZ=[\z_1, \z_2,\cdots,\z_{k_1}]$. Also, let $\Y_{\text{ext}}=[\Y,\y_{k_1+1},\cdots,$ $\y_{k_2}]$ be an $\text{OA}(N,k_2,2,t)$ extension of $\Y$.
By OD-equivalence, $\Y$ is isomorphic to $\tilde{\ZZ}=[\z_j\odot\z_1,
 \z_j\odot\z_2, \cdots, \z_j\odot\z_{j-1}, \z_j, \z_j\odot\z_{j+1}, \cdots,  \z_j\odot\z_{k_1}] $
for some $j \in \{1,\ldots,k_1\}$.
Since $\tilde{\ZZ}$ and $\Y$ are isomorphic, $\tilde{\ZZ}$ is an $\text{OA}(N,k_1,2,t)$ and extends to an $\text{OA}(N,k_2,2,t)$,
$\tilde{\ZZ}_{\text{ext}}=[\z_j\odot\z_1, \z_j\odot\z_2, \cdots, \z_j\odot\z_{j-1}, \z_j, \z_j\odot\z_{j+1}, \cdots,  \z_j\odot\z_{k_1}, \z_{k_1+1},\cdots,\z_{k_2}].$
Let $\tilde{\tilde{\ZZ}}_{\text{ext}}=[\z_1, \z_2, \cdots, \z_{j-1}, \z_j,  \z_{j+1},  \cdots,   \z_{k_1}, \z_j\odot\z_{k_1+1}, \cdots, \z_j\odot\z_{k_2}]$.
It follows that $\tilde{\tilde{\ZZ}}_{\text{ext}}$ is OD-equivalent to $\tilde{\ZZ}_{\text{ext}}$. Furthermore, each $J_r(\ell)$ value of $\ZZ$ and $\tilde{\tilde{\ZZ}}_{\text{ext}}$ equals a  $J_r(\ell)$ or $J_{r-1}(\ell\backslash \{j\})$ value of $\tilde{\ZZ}$ and $\tilde{\ZZ}_{\text{ext}}$ if $r$ is even and a $J_r(\ell)$ or $J_{r+1}(\ell\cup\{j\})$ value of $\tilde{\ZZ}$ and $\tilde{\ZZ}_{\text{ext}}$
if $r$ is odd. This implies that $$ J_r(\ell)\left(\tilde{\tilde{\ZZ}}_{\text{ext}} \right)=\begin{cases}0\quad \mbox{ for } r\leq t \mbox{ if } t \mbox{ is even,}\\
0 \quad \mbox{ for } r\leq t-1 \mbox{ if } t \mbox{ is odd},
\end{cases}
$$ where $ J_r(\ell)\left(\tilde{\tilde{\ZZ}}_{\text{ext}}\right)$ are the $J_r(\ell)$ values  of $\tilde{\tilde{\ZZ}}_{\text{ext}}$.
  Hence,  $\tilde{\tilde{\ZZ}}_{\text{ext}}$ is the desired extension, and $\ZZ=[\z_1, \z_2,  \cdots, \z_{k_1}]$ has the desired strength. Example~\ref{ex:ex} shows that when $t$ is odd $\ZZ$ does not have to be an OA$(N,k_1,2,t)$.
\qed
\end{pf}


Let $M$ be a set of $\oant$.
Next, we show how to remove OD-equivalent $\text{OA}(N,k,2,$ $t)$ from $M$ 
by solving an equivalent graph isomorphism problem~\cite{McKay}. 
First, a column of $1$s is concatenated to each $\text{OA}(N,k,2,t)$ in $M$. Second, each resulting design is converted to the graph defined by McKay~\cite{McKay}. Third, the set of resulting graphs is reduced to a subset of all non-isomorphic graphs using the \texttt{shortg} utility in \texttt{nauty}~\cite{McKay13}. Finally, the resulting non-isomorphic graphs are converted to a set of all non-OD-equivalent OAs. We call this process the {\em OD-equivalence reduction} of OAs. We also call the process of removing isomorphic OAs  from a set $M$ of OAs to obtain a set of all non-isomorphic OAs in $M$ the {\em isomorphism reduction} of OAs.

Given a set $M$  of all non-OD-equivalent OAs, a set of all non-isomorphic $\oant$ can be extracted from $M$ by following Steps 1-4 below.
\begin{enumerate}
\item For each $\Y \in M$, compute $\y_i \odot [\1,\Y]=[\y_i,\y_i \odot \y_1,\cdots,\y_i \odot \y_k]$ for $i=1,\ldots,k$, where $\y_i$ is the $i$th column of $\Y$.
\item Convert each resulting design as well as each $[\1, \Y$] for $\Y \in M$ to the graph $G_1$ for $N$ row, $k$ column designs in Ryan and Bulutoglu~\cite{Ryan} after deleting its column of all $1$s.
\item Reduce this set of graphs to a subset of non-isomorphic graphs with \texttt{shortg}~\cite{McKay}.
\item Convert each graph back to its corresponding design~\cite{Ryan}.
\end{enumerate}
We do not report the time needed to extract a set of all non-isomorphic $\oant$ from a set of non-OD-equivalent $\oant$ in Table~\ref{tab:enum} since it is insignificant compared to the time needed to obtain a set of all non-OD-equivalent $\oant$.

\section{Symmetric ILP and  B\&B with isomorphism pruning}\label{sec:isop}
 \setcounter{dummyt}{0} 
\setcounter{dummym}{0} 
\setcounter{dummya}{0} 
\setcounter{dummyc}{0}
\setcounter{dummyl}{0}
\setcounter{dummyd}{0}
\setcounter{dummye}{0}
We first define the symmetry group of an ILP and the notion of a symmetric ILP.
For some integer $p \geq 1$,  matrices $\A$ and $\B$ with $n$ columns, column vectors $\bb$ and $\dd$ with the same number of rows as $\A$ and $\B$, and $\cc \in \mathbb{R}^n$, let  $\mathcal{F}$ be the feasible set of ILP
\begin{equation}\label{eqn:geneqILP}
\begin{array}{ll}
&\mbox{min}  \ \cc^{\top}\x \quad \mbox{subject to:}\\
 &\A \x=\bb, \ \B \x\leq\dd,\\
 &\x \in \{0,\ldots,p\}^n ;
\end{array}
\end{equation}
let $S_n$ be the set of all permutations of the indices $\{1,\ldots,n\}$; and let
\begin{equation*}
G=\{\pi \in S_n \  | \ \pi(\x) \in \mathcal{F} \   \text{and }\cc^{\top}\pi(\x)=\cc^{\top}\x\ \text{ for all } \x \in \mathcal{F}\},
\end{equation*} where
$$\pi\left((x_1,x_2,\ldots, x_n)^{\top}\right)= (x_{\pi^{-1}(1)},x_{\pi^{-1}(2)},\ldots,x_{\pi^{-1}(n)})^{\top} \quad \forall \+ \x \in \mathbb{R}^n.$$
Then the group $G$ is called the {\em symmetry group} of  ILP~(\ref{eqn:geneqILP})~\cite{Margot2010}.
If $|G|>1$, then ILP~(\ref{eqn:geneqILP}) is called {\em symmetric}, and each non-identity element in $G$ is called a {\em symmetry} of ILP~(\ref{eqn:geneqILP})~\cite{Margot2010}. 
Let $H$ be a subgroup of   the symmetry group of an ILP.  Then
 two solutions  $\x_1$ and $\x_2$ of the ILP are called {\em isomorphic under the action of $H$} if there exists an $h\in H$   such that $\x_1=h(\x_2)$. 

Bulutoglu and Margot~\cite{BulutogluM} showed that finding all $\oan$ is equivalent to finding all non-negative integer solutions to a symmetric ILP.   A basic B\&B based on linear programming (LP)  is a standard backtrack algorithm for finding a solution to an ILP, where partial solutions (nodes) in its search tree are pruned by solving LP relaxations, see Chapter 7 of Wolsey~\cite{Wolsey}. When each feasible solution of an ILP is also optimum,  a basic B\&B algorithm that prunes a node of its search tree if and only if the LP relaxation at that node is proven to be infeasible by  the dual simplex algorithm (see~\cite{Banciu}  for a treatment of the dual simplex algorithm) can  be used to find all feasible solutions. We take this approach  as the objective function of each ILP considered in this paper is the zero function, so each solution is optimum. 
 However, the presence of symmetry in an ILP causes  a basic
  B\&B  to visit many nodes (subproblems) with the same LP relaxation objective function value and feasibility status.  This requires improving a basic B\&B algorithm to decrease the number of  such subproblems. Such an improvement was developed by Margot~\cite{Margot02,Margot03a,Margot03b,Margot07} and used in~\cite{INFORMS} to improve the lower bound from $65$ to $71$ for the $6$ matches instance of the famous football pool problem.
Next, we describe this algorithm by using the same notation as in~\cite{INFORMS} adapted for ILP~(\ref{eqn:geneqILP}). 

Let each solution to ILP~(\ref{eqn:geneqILP}) be also optimum.
Let $\tau$   be the basic B\&B tree after implementing a basic  B\&B to find all (optimum) solutions to ILP~(\ref{eqn:geneqILP}). (The search tree of a basic B\&B algorithm depends on the branching strategy used~\cite{BALA}).
Let $j\in \{0,\ldots,p\}$, and
$$F_j^a=\{i\ | \ x_i \text{ has been fixed to } j \text{ by branching leading up to node $a$}\}.$$ 
Nodes $a$ and $b$ represent the same node if 
$F_j^a=F_j^b$ for $j=0,\ldots, p$. Let $H$ be a subgroup of the symmetry group $G$ of ILP~(\ref{eqn:geneqILP}).
Nodes $a$ and $b$ are {\em isomorphic} under the action of $H\leq G$ if 
$$\exists \+ \pi \in H \text{ such that } \pi (F_j^a)=F_j^b \text{ for } j=0,\ldots,p,$$ where  $\pi(S)=\{\pi(i)\ |\ i \in S\}$ 
for $S \subseteq \{1,\ldots,n\}$.
The action of $H$ on the nodes of $\tau$ partitions the set of all possible nodes into isomorphism classes.
Let a basic B\&B algorithm on ILP~(\ref{eqn:geneqILP}) be implemented in such a way that   exactly one class representative of each node  in $\tau$ under the action of $H$ is kept. Then removing isomorphic nodes in this way is called  {\em isomorphism pruning} and the resulting algorithm is called {\em B\&B with isomorphism pruning}~\cite{INFORMS,Margot02,Margot03a,Margot03b,Margot07}. The following lemma shows that a B\&B with isomorphism pruning algorithm can be used to find all non-isomorphic solutions to an ILP with the zero objective function under the action of $H$.
\begin{lem}
Let a B\&B  with isomorphism pruning be implemented on ILP~(\ref{eqn:geneqILP}) with $\cc=\zz$ using the group $H \leq G$, where $\zz$ is the all zeros vector, and $\tau'$  be the resulting B\&B tree after  isomorphism pruning. Then the resulting set of all leaves (feasible and optimal solutions) in $\tau'$ constitute a set of all non-isomorphic solutions under $H$.
\end{lem}
\begin{pf}
Let $\tau$ be the resulting basic B\&B tree  if isomorphism pruning is not implemented. 
Since  $\cc=\zz$, every feasible solution is optimum and each leaf of $\tau'$ and $\tau$ is a solution. The result now follows from the way $\tau'$ is defined based on $\tau$. \qed
\end{pf}
Randomly picking class representatives at each depth in $\tau$ does not always result in a tree where the resulting graph may be disconnected.
Hence,  for a valid B\&B with isomorphism pruning algorithm the class representatives of nodes  should be picked in such a way that the resulting graph after deleting the non-class representative nodes is a tree.
On the other hand, if at each node the minimum index non-fixed variable is selected for branching (called the {\em minimum index branching}), removing each node that is not lexicographically minimum within its class results in a tree and hence a valid algorithm~\cite{Margot07}.
This follows from the fact that if $\tau$ is constructed by minimum index branching, then each unique lexicographically minimum node in its class in $\tau$ under the action of $H$ has a unique lexicographically minimum node as a parent under the same action of $H$. 

Margot~\cite{Margot07}  developed a group-theory-based algorithm for deciding  whether a given node in $\tau$ is lexicographically minimum within its class. This algorithm is used within the B\&B  with isomorphism pruning algorithm in~\cite{Margot07} to prune nodes that are not lexicographically minimum within their isomorphism classes without resolving LP relaxations.
On average, using larger subgroups $H$ for isomorphism pruning results in faster solution times. This is because after each branching decision isomorphism pruning prunes a
larger number of isomorphic nodes with a larger subgroup without resolving  an LP relaxation, and on average deciding whether a node is lexicographically minimum in its class by using the 
group-theory-based algorithm in~\cite{Margot07}
is faster than  resolving.
 However, there are exceptions. In particular, a group that has the form 
$S_{r_1}\times S_{r_2}\times \cdots \times S_{r_h}$  with a coordinate action as defined in Pfetsch and Rehn~\cite{REHN}, where several $r_i\geq 3$, is sometimes a worst-case example for group algorithms that determine whether a given node in the B\&B tree is lexicographically minimum under the action of the group~\cite{Margot07,REHN}. Hence, such symmetries should be removed by using inequalities, and the group-theory-based algorithms should be used to exploit the remaining symmetries.
 \section{ILP formulations for classifying {\boldmath $\oant$}}\label{sec:formulations}
 \setcounter{dummyt}{0} 
\setcounter{dummym}{0} 
\setcounter{dummya}{0} 
\setcounter{dummyc}{0}
\setcounter{dummyl}{0}
\setcounter{dummyd}{0}
\setcounter{dummye}{0}
Bulutoglu and Margot~\cite{BulutogluM} used the B\&B with isomorphism pruning algorithm in~\cite{Margot07} for finding all non-isomorphic solutions to one ILP in $s^k$ variables to directly classify $\oan$. The full ILP formulation  from~\cite{BulutogluM} is defined by $\binom{k}{t}s^t$ equality constraints and an unknown solution vector $\x:=(x_1,x_2,\ldots,x_{s^k})^{\top}$ whose $\left (\sum_{j=1}^k i_j s^{k-j}+1\right)$th entry is the number of times the  symbol combination $(i_1,i_2,\ldots,i_k) \in \{0,\ldots,s-1\}^k$ appears in the rows of the sought after $\oan$. Next, we provide an improved full ILP formulation  whose LP relaxation does not have redundant constraints. The improvement stems from having $\binom{k}{t}s^t -\left(\sum_{q=0}^t(s-1)^q\binom{k}{q} \right)$ fewer  equality constraints, yet the new set of constraints are equivalent to those in~\cite{BulutogluM}. This improved full ILP formulation is
 
\begin{align}
 &\min  \ 0 \quad \text{subject to:}\hspace{0cm} \nonumber\\
  & \sum_{\{i_{1},i_2,\ldots,i_{k}\} \backslash \{i_{j_1}, i_{j_2},\ldots,i_{j_q}\} \in \{0,\ldots,s-1\}^{k-q} } {x_{[i_1 s^{k-1}+i_2s^{k-2}+\ldots+i_k s^{k-k}+1]}}= \frac{N}{s^q},\label{ilp:BF}\\
 & 1 \leq x_1, \  \ \x \in \{0,\ldots, p_{\max}\}^{s^k} \nonumber
\end{align}
for   $q=0,\ldots, t$    and  each vector of  $q$ indices  $(i_{j_1}, i_{j_2},\ldots,i_{j_q})\in \{0,\ldots,s-2\}^{q}$ for each $\{j_1,j_2,\ldots,j_q\} \subseteq \{1,\ldots,k\}$, where
 $p_{\max}$ is computed by solving the ILP in Lemma
5 from Bulutoglu and Margot~\cite{BulutogluM}.

The following lemma 
(which follows directly from the definition of an $\oan$) and definition are used to formulate the ILP feasibility problem resulting from extending an $\oankm$ to $\oan$ without losing isomorphism classes. 
\begin{lem}\label{lem:Y'}
Let $\Y$ be an $N$ row, $k$ column, $s$-symbol  array with
columns $\{\y_1,\y_2,\ldots,$ $\y_k\}$. Let $\Y'$ be the $N \times
(s-1)k$ matrix with columns $\{\y'_1,\y'_2,\ldots,\y'_{(s-1)k}\}$, where for
$j=1, \ldots, k$ and $r=1, \ldots, s-1$ the $i$th entry of
$\y'_{(s-1)(j-1)+r}$ is $1$ if the $i$th entry of $\y_j$ is $r-1$ and 
$0$ otherwise.  Then $\Y$ is an $\oan$ if and only if
\begin{equation}\label{eqn:jy'01}
\sum_{i=1}^Ny'_{ih_1}y'_{ih_2}\ldots\,y'_{ih_q}=\frac{N}{s^q} \\
\end{equation}
for any $q$ columns $\{\y'_{h_1},\y'_{h_2},\ldots,\y'_{h_q}\}$ of
$\Y'$ such that $\lceil h_{i'}/(s-1) \rceil \ne \lceil h_{j'}/(s-1) \rceil$ for all $q=1,\ldots,t$
and $1 \leq i'<j' \leq q$.
\end{lem}
\begin{defn}\label{def:cor}
Let $\Y'$ from Lemma~\ref{lem:Y'} be a solution to the system of $\sum_{q=1}^t (s-1)^{q}\binom{k}{q}$ equations~(\ref{eqn:jy'01}). Also, let $\Y$ be the $\oan$ obtained from $\Y'$ with columns
\begin{equation*}
\y_{j}=\sum_{r=1}^{s-1}{(r-1)\y'_{(s-1)(j-1)+r}}+(s-1)\left(\1-\sum_{r=1}^{s-1}{\y'_{(s-1)(j-1)+r}}\right)\quad
\end{equation*}
for each $j=1,\ldots,k$. Then $\Y$ is called {\em the $\oan$ corresponding to $\Y'$}.
\end{defn}

 Given an input OA$(N,k-1,s,t)$ say $\Y$ that is being extended to $\oan$, construct $\Y'$ from $\Y$ using Lemma~\ref{lem:Y'}, and add $s-1$ columns of binary variables $\xx_r=(x_{1r},x_{2r},\ldots,x_{Nr})^{\top}$ for $r=1,\ldots,s-1$ to $\Y'$. Then an  ILP
 in $N(s-1)$ binary variables  for solving this extension problem is 
\begin{align}\label{ilp:oa}
& \min  \ 0 \quad \text{subject to:} \nonumber\\
&  \sum_{i=1}^Nx_{ir}=\frac{N}{s}, \nonumber \\
&\sum_{i=1}^Ny'_{ih_1}y'_{ih_2}\ldots\,y'_{ih_{q-1}}x_{ir} =\frac{N}{s^q},\\
&x_{1,1}=1, \ \ \ \sum_{r=1}^{s-1}x_{ir} \leq 1, \ \ \ \sum_{m=1}^{r}{(x_{i''m}-x_{(i''+1)m})}
\ge 0, \nonumber\\
 &  \xx_r \in \{0,1\}^{N}, \quad i=1,\ldots,N, \quad r=1,\ldots,s-1, \quad q=2,\ldots,t \nonumber
\end{align}
 for each pair of identical rows whose indices are of the form $( i'',i''+1)$  in the input OA $\Y$, where $1 \leq i'' <  i''+1 \leq N$
 and for any $q-1$ columns $\y'_{h_1}, \y'_{h_2},\ldots,\y'_{h_{q-1}}$
as in Lemma~\ref{lem:Y'} with the last $s-1$ columns of $\Y'$  deleted.  Each solution matrix $[\xx_1,\x_2,\cdots,\x_{s-1}]$ to
ILP~(\ref{ilp:oa}) has a corresponding $\oan$. 

The $\sum_{q=1}^t (s-1)^{q}\binom{k-1}{q-1}$ constraints in ILP~(\ref{ilp:oa}) for $q=1,\ldots,t$ with $N/s^q$ on the right hand side imply that the entries in any $t$ subset of  columns including the new column will be such that any $t$-tuple in  $\{0,\ldots,s-1\}^t$  will appear $\lambda$ times as rows of the $t$ subset of columns. Fixing the first entry of the new column at symbol $0$  by setting $x_{1,1}=1$ and ordering column symbols in the new column over identical rows within the input OA by requiring $\sum_{m=1}^{r}({x_{i''m}-x_{(i''+1)m}})\geq0$ for $r=1,\ldots,s-1$ boosts speed by removing some of the symmetries  without losing isomorphism classes. 

Next, a hybrid  version of ILPs~(\ref{ilp:BF}) and~(\ref{ilp:oa})   is defined. It reduces the number of variables in the improved full ILP formulation  from $s^k$ to $hs$, where $h \le N$ is the number of distinct rows in the input $\oankm$ that is being extended. Let the rows of an input OA $\Y$ be ordered lexicographically so that identical rows appear next to each other. Let $r_i$ for $i=1,\ldots,s^{k-1}$ be the number of times the $i$th row of $\D_s^{k-1}$ appears in the input OA, where the set of rows of $\D_s^{k-1}$ is  $\{0,\ldots,s-1\}^{k-1}$ and the rows of  $\D_s^{k-1}$ are ordered lexicographically. Let $1=i_1 < i_2 < \cdots < i_h \leq s^{k-1}$ be such that $r_{i_{\ell}}>0$ for $\ell=1,\ldots,h$.

Construct $\Y'$ from $\Y$ as in  Lemma~\ref{lem:Y'}. Let $1=p_1 < p_2 < \cdots < p_h \leq N$ be such that the $p_{\ell+1}$th row in the input OA is the first row different from the $p_{\ell}$th row. Let $\x$ be the indicator vector of the sought after $\oan$ as in ILP~(\ref{ilp:BF}) and $\x^{h(s-1)}$ be the resulting vector with $h(s-1)$ entries after deleting all the entries of $\x$ that do not appear in the constraints of the following ILP 
\\
\begin{align}\label{ilp:oahybrid}
&  \min  0 \quad \text{subject to:}\nonumber\\
&   \sum_{\ell=1}^h x_{(i_{\ell}-1)s+j} =\frac{N}{s},\notag \\
&\sum_{\ell=1}^hy'_{p_{\ell}\alpha_1}y'_{p_{\ell}\alpha_2}\ldots\,y'_{p_{\ell}\alpha_{q-1}}x_{(i_{\ell}-1)s+j} =\frac{N}{s^q},\\
&1 \leq x_{1}, \ \ \ \sum_{j=1}^{s-1}x_{(i_{\ell}-1)s+j}\leq r_{i_{\ell}}, \nonumber\\
 &   x_{(i_{\ell}-1)s+j} \in \{0,\ldots,\min(r_{i_{\ell}},p_{max})\},   \nonumber\\
&\text{for }\ell=1,\ldots,h, \ j=1,\ldots,s-1, \ q=2,\ldots,t, \nonumber
\end{align}
and for any $q-1$ columns  $\{\y'_{\alpha_1},\y'_{\alpha_2},\ldots,\y'_{\alpha_{q-1}}\}$
as in Lemma~\ref{lem:Y'} with the last $s-1$ columns of $\Y'$  deleted. Each solution $\x^{h(s-1)}$ to
ILP~(\ref{ilp:oahybrid}) has a corresponding $\oan$.

One fundamental difference between ILPs~(\ref{ilp:oa}) 
and~(\ref{ilp:oahybrid}) is that ILP~(\ref{ilp:oa}) depends on the ordering of the rows in the input OA and ILP~(\ref{ilp:oahybrid}) does not.
It is well known that different ILP formulations of the same ILP problem can result in orders of magnitude variability in solution times~\cite{BALA}, and this is especially the case for ILPs with large symmetry groups when using a symmetry exploiting algorithm~\cite{Margot2010}. Consequently,  
our two ILP formulations~(\ref{ilp:oa}) and~(\ref{ilp:oahybrid}) were needed because neither was sufficient to extend all non-isomorphic OA$(N,k-1,2,t)$ to  all non-isomorphic $\oant$ for all $k$ such that an OA$(N,k-1,2,t)$ exists.
\section{ILP-based methods for classifying {\boldmath  $\oan$}}\label{sec:enum}
Speed improvements gleaned from using larger groups  with the isomorphism pruning  algorithm of Margot~\cite{Margot07} is the main motivation to determine the largest  groups that are viable to use in enumerating all non-isomorphic solutions of ILP~(\ref{ilp:oa}) and ILP~(\ref{ilp:oahybrid}) to extend an  OA$(N,k-1,s,t)$ to $\oan$ without losing isomorphism classes. Such  groups are found and used in methods ILP-extension and Hybrid.
 

For an $\oan$ $\Y$, let $\Y''$ be the $N\times sk$ matrix obtained from $\Y'$ in Lemma~\ref{lem:Y'} by concatenating columns $\1_N-\sum_{r=1}^{s-1}{\y'_{(s-1)(j-1)+r}}$ for $j=1,\ldots, k$. Let $$G_s^k=\{\pi \in S_N\ | \mbox{ there exists } \sigma
\in S_{sk} \mbox{ such that } \Y''(\pi,\sigma)=\Y''\},$$ where  $\Y''(\pi,\sigma)$
is the resulting matrix when the rows and columns of $\Y''$ are
permuted according to $\pi$ and $\sigma$, respectively. Let $H_s^{k}$ be the maximum subgroup of
$G_s^{k}$ that does not have a non-identity element that sends a row's index of $\Y$ to the index of an identical row. Next, we describe the ILP-extension method.

\begin{meth}[ILP-Extension] \label{alg:BR}
Input: An OA$(N,k-1,s,t)$ $\Y$. 

 Find all non-isomorphic solutions to ILP~(\ref{ilp:oa})
under the action of the group $H_s^{k-1}$ using the B\&B with isomorphism pruning algorithm~\cite{Margot07}.

 Output: The set $M$ of all such $\oan$.
\end{meth}

\begin{thm}\label{thm:BR3idd}
The output $M$ from Method~\ref{alg:BR} contains at least one
  $\oan$ from each isomorphism class of $\oan$ that extends the input $\Y$.
 \end{thm}
\begin{pf}
Each element of $M$ from Method~\ref{alg:BR} is an $\oan$ by Lemma~\ref{lem:Y'}.
Let Method~\ref{alg:BR} be implemented with an $\oankm$ input $\Y_1$. We show that the set $M$ contains at least one copy of each non-isomorphic extension of $\Y_1$ to an $\oan$. For each $g \in H_{s}^{k-1}$, $g(\Y'_1)$ corresponds to an $\oankm$.  Let $\Y^\sigma_1$ denote this $\oankm$. Then $\Y^\sigma_1$ is isomorphic to $\Y_1$ with an isomorphism $\sigma$ that does not involve a permutation of the rows (i.e., $\Y^\sigma_1$ is obtained from $\Y_1$ by applying an element  $\sigma \in S_s\wr S_k$ to the columns of $\Y_1$). Now, each solution $\xx$ of  ILP~(\ref{ilp:oa}) corresponds to an extension of $\Y_1$ that is an $\oan$. If $[\Y_1, \y]$ is this extension, then $g(\xx)$ also solves this ILP and corresponds to $[\Y^\sigma_1, \y]$, so $[\Y_1, \y]$ and $[\Y^\sigma_1, \y]$ are isomorphic.
\qed
\end{pf} 
Let the set of rows of $\D_s^k$ be $\{0,\ldots,s-1\}^k$ and the rows of 
$\D_s^k$ be ordered lexicographically. Let $G_{s,k}$ be the group of  
 all permutations of columns of $\D_s^k$ followed by (possibly different)  permutations of $\{0,\ldots,s-1\}$ within each column. Then $G_{s,k}$ acts on the rows of $\D_s^k$ and consequently acts on the variables of ILP~(\ref{ilp:BF}). In fact, $G_{s,k}\cong S_s \wr S_k$ is the paratopism group that acts on $\oan$~\cite{Egan}. Hence, each element of $G_{s,k}$ maps a solution $\x$ of ILP~(\ref{ilp:BF})
corresponding to an $\oan$ to that of an isomorphic OA; this
implies that the B\&B with isomorphism pruning algorithm using $G_{s,k}$ 
returns the unique $\oan$ that is lexicographically minimum in
rows within the orbit of $G_{s,k}$ per isomorphism class, where for a group $G$ and  $\vv \in \mathbb{R}^n$,  $\{\pi(\vv)\ |\ \pi \in G\}$ is the  {\em orbit} of $\vv$ under the action of $G$. The constraint $1\leq x_1$ is valid for the sought after isomorphism class representatives, and its
inclusion boosts speed.
Hence, the following method finds a set of all non-isomorphic $\oan$.
\begin{meth}[Full]\label{alg:BF} Input: $N$, $k$, $s$,
$t$.

 Find all non-isomorphic solutions to ILP~(\ref{ilp:BF}) under the action of the group $G_{s,k}$ using the
  B\&B with isomorphism pruning algorithm~\cite{Margot07}.

Output: A set $M$ of all non-isomorphic $\oan$. (When $s=2$, optionally reduce $M$ to retain one representative from each OD-equivalence class.)
\end{meth}

Let $H_{s,k}$ be the maximum subgroup of $G_{s,k}$ that maps the input $\oankm$ to itself. Next, we describe the method Hybrid.

\begin{meth}[Hybrid] \label{alg:hybrid}
Input: An OA$(N,k-1,s,t)$ $\Y$.

 Find all non-isomorphic solutions to ILP~(\ref{ilp:oahybrid}) under the action of the group $H_{s,k}$ using the
 B\&B with isomorphism pruning algorithm~\cite{Margot07}.

 Output: The set $M$  of all such $\oan$.
\end{meth}

\begin{thm}\label{thm:hybrid4}
An $\oan$ from each isomorphism class of $\oan$ that extends $\Y$ appears at least once in the output of  Method~\ref{alg:hybrid}.
\end{thm}
\begin{pf}
The indicator variables for  ILP~(\ref{ilp:BF}) with $s^k$ variables for finding
all $\oan$ are indexed with the elements of $\{0,\ldots,s-1\}^k$. The
action of $G_{s,k}$ on the elements of  $\{0,\ldots,s-1\}^k$ defines the action
of $G_{s,k}$ on the indicator variables. All elements in $G_{s,k}$
are viable to use in isomorphism pruning, except those that produce
a different $\oankm$ input OA when $G_{s,k}$'s action is
restricted to the first $k-1$ columns. This is true since the
elements that are not viable send a feasible solution (i.e., a
feasible extension of the input OA) to a potentially infeasible point
(i.e., an infeasible extension  whose first $k-1$ columns up to a permutation of rows is not the input
design). The elements that are viable to use in isomorphism pruning  
send a row that appears $r_i$ times in the input OA to a row
that appears the same number of times. These elements form the group
$H_{s,k}$.
\qed
\end{pf}
The solution times for ILP~(\ref{ilp:BF}) as well as those for the ILPs~(\ref{ilp:oahybrid})  do depend on the ordering of the variables. Different orderings of variables result in drastically different solution times. This is mainly because changing the order of the variables changes the order in which branching is implemented. (With the B\&B with isomorphism pruning algorithm, we always opted to branch in the index order of the variables.) The lexicographical ordering of the variables given in ILPs~(\ref{ilp:BF}) and~(\ref{ilp:oahybrid}) resulted in substantially faster solution times when  compared to random orderings.

Let $f(N,s,t)$ be the largest $k$ for which an $\oan$ exists~\cite{Hedayat}. For $k=t,\ldots,f(N,s,$ $t)$, let $n_k$ be the number of
non-isomorphic $\oan$ and $T_k=\{\Y_1, \Y_2, \ldots, \Y_{n_k}\}$ be a
set of non-isomorphic $\oan$. Alternatively, $T_k$ can be a set of non-OD-equivalent OA$(N,k,s,$ $t)$ when $s=2$.  In what follows, we describe a generic method (Basic Extension) for obtaining $T_k$ from $T_{k-1}$.
\begin{meth}[Basic Extension] \label{alg:basic} Input: $N$, $k$, $s$, $t$,
and $T_{k-1}$.  Set $\ell:=1$.
\begin{enumerate}

\item Using ILP-extension or Hybrid obtain a set $M_{\ell}$ of $\oan$ such that the first $k-1$ columns
are the same as those of $\Y_{\ell} \in T_{k-1}$ and at least one representative
from each isomorphism class (OD-equivalence class when $s=2$) of such $\oan$ is
included.\label{item:extend}

\item Increment $\ell:=\ell+1$ and then repeat Step~\ref{item:extend} if $\ell \leq
n_{k-1}$.\label{item:iterate}

\item \label{item:M} Set $M:=\bigcup_{\ell=1}^{n_{k-1}}M_{\ell}$. Reduce $M$ to $T_k$ with one representative from each isomorphism class (OD-equivalence class when $s=2$) of
$\oan$ in $M$.

 Output: $T_k$.\label{item:nauty}
\end{enumerate}
\end{meth}
  Let $T_t$ consisting of one OA$(N,t,s,t)$ be initialized so that its rows are ordered lexicographically.  Given $T_{k-1}$, $T_{k+m}$ is obtained after applying the Basic Extension method  $m+1$ times. Alternatively, if $k>t+1$, Bulutoglu and Margot~\cite{BulutogluM} (with OD-equivalence reduction when $s=2$)  can be used to directly obtain input $T_{k-1}$ to the Basic Extension method. Our implementations of Basic Extension use the graph-based method involving the program \texttt{nauty}~\cite{McKayP,McKay13} as described in Ryan and Bulutoglu~\cite{Ryan} for isomorphism reductions in Step~\ref{item:nauty}. 
When $s=2$ OD-equivalence reductions are implemented by using the graph-based  method  in Section~\ref{sec:fold}.

Implementing OD-equivalence in place of OA-isomorphism reductions greatly reduces the number of ILPs in Step~\ref{item:extend} of Method~\ref{alg:basic} when $t$ is even; see Table~\ref{tab:enum}. On the other hand, OD-equivalence reductions are not worthwhile when $t$ is odd because of Theorem~\ref{thmOD}(ii). First, a strength $t$ array satisfies a  strict subset of the constraints that a strength $t+1$ array satisfies, so there often tends to be many more $\oant$ that are not of strength $t+1$ than $\text{OA}(N,k,2,t+1)$. In addition, an $\oant$ with odd $t$ that is not of strength $t+1$ can only be OD-equivalent but not isomorphic to an $\text{OA}(N,k,2,t-1)$ that is not of strength $t$, but such $\text{OA}(N,k,2,t-1)$ are already excluded. Hence, the only possible OD-equivalence reductions involve the much smaller number of $\oant$ with strengths exceeding $t$ when $t$ is odd.

In ILP~(\ref{ilp:oa}), the constraints 
\begin{equation}\label{eqn:symbreak}
\sum_{m=1}^{r}({x_{i''m}-x_{(i''+1)m}})\geq0  \quad \text{ for } r=1,\ldots,s-1
\end{equation}
remove all the symmetry due to identical rows provided that each distinct row $\z$ in the input OA appears in rows
$i_{\z},i_{\z}+1,\ldots, i_{\z} +r_{\z}-1$ for some $i_{\z} \in \{1,\ldots,  N\}$, where $r_{\z}$  is the number of times $\z$ appears. 
This property is indeed satisfied by all the OAs produced by Full,  followed by Hybrid,  
or by Full,  followed by Hybrid,  followed by ILP-extension, and in our classifications,  all the input OAs to ILP-extension are produced in one of these two orders.  

Removing some of the symmetry using constraints~(\ref{eqn:symbreak}) is necessary. Otherwise,
 the B\&B with isomorphism pruning algorithm~\cite{Margot07} stalls. 
In fact, for the $s=2$ cases, the group of all symmetries removed by introducing constraints~(\ref{eqn:symbreak})
is isomorphic to $S_{r_1}\times S_{r_2}\times \cdots \times S_{r_h}$, 
where $r_i$ is the number of times the $i$th distinct row is repeated in the input OA.  As discussed in Section~\ref{sec:isop}, it is desirable to remove  these symmetries using inequalities and use group-theory-based algorithms to exploit the remaining symmetries.

For each of the ILPs~(\ref{ilp:BF}),~(\ref{ilp:oa}), and~(\ref{ilp:oahybrid}) the objective function $0$ is picked since our goal is to find all non-isomorphic solutions to its system of constraints not only those that minimize a particular objective function. We could have also picked the left hand side of one of its equality constraints that also does not decrease the size of its group of all exploitable symmetries.
 However, using $0$ instead of other viable objective functions is faster. Next, we explain why this happens.
 
The   B\&B with isomorphism pruning algorithm~\cite{Margot07} and  an implementation of a basic B\&B in {\tt CPLEX12.5.1} both use the {\tt CPLEX12.5.1} implementation of the dual simplex algorithm to resolve the LP relaxations at the nodes in their search trees.
When $0$ is the objective function of an ILP then all of its reduced costs are zero. Hence, all the ratios in the ratio test of the dual simplex algorithm are equal to zero and the dual of the LP relaxation is degenerate. This causes the {\tt CPLEX12.5.1} implementation
to perturb the coefficients in the objective function and the perturbation speeds up the LP resolve. When another viable objective function is used {\tt CPLEX12.5.1} does not necessarily apply such a perturbation.
However, such a perturbation would also be beneficial when other viable 
objective functions are used.  This is mainly because the variability of the coefficients in these viable objective functions is small, and when this is the case the dual simplex algorithm tends to be slower~\cite{cplex}.

It might be argued that the only difference between the MCS algorithm and the ILP-based methods ILP-extension or Hybrid is that
the MCS algorithm does not solve LP relaxations in its search tree. However, 
there are other major differences. First, when ILP-extension or Hybrid is used, isomorphism (OD-equivalence) reductions are achieved with \texttt{nauty} as opposed to the lexicographically minimum in columns (LMC) check algorithm~\cite{Schoen}. Second, 
 when one of ILP-extension or Hybrid is implemented with its isomorphism pruning, it uses both constraints and a 
 group-theory-based algorithm for rejecting isomorphic partial solutions (nodes in the B\&B  tree), 
   whereas the MCS algorithm uses only constraints for the same purpose.

\section{Computational results and performance comparisons}\label{sec:results}
 \setcounter{dummyt}{0} 
\setcounter{dummym}{0} 
\setcounter{dummya}{0} 
\setcounter{dummyc}{0}
\setcounter{dummyl}{0}
\setcounter{dummyd}{0}
\setcounter{dummye}{0}
The generic method Basic Extension in Section~\ref{sec:enum} was used with different methods to extend the known $\oant$ catalog.  Classification results are summarized in Table~\ref{tab:enum}, including the number of OD-classes whenever OD-equivalence reductions were used. In such cases, Steps 1-4 from Section~\ref{sec:fold} were necessary to obtain a set of all non-isomorphic OAs from a set of all non-OD-equivalent ODs. The new results are a set of all non-isomorphic $\text{OA}(160,9,2,4)$, no $\text{OA}(160,10,2,4)$ exists, a set of all non-isomorphic $\text{OA}(176,k,2,4)$ for $k=5,6,7,8,9$, and no $\text{OA}(176,10,2,4)$ exists, so $f(160,2,4)=f(176,2,4)=9$ is now established. These were the smallest unknown $2$-symbol strength $4$  cases. Schoen, Eendebak, and Nguyen~\cite{Schoen} did classify $\text{OA}(160,k,2,4)$ for $k=5,6,7,8$, (the largest $2$-symbol, strength $4$ cases classified so far), and their  results also match Table~\ref{tab:enum}.  

 Isomorphism pruning  was always turned off with the ILP-extension method. This is mainly because most ILPs of this method  had little or no exploitable symmetry, and the time spent in finding and exploiting these symmetries was not  worth the effort. For  the other  methods, the time spent in finding and exploiting these symmetries paid off well, and  isomorphism pruning was turned off only if the size of the exploitable  symmetry group was strictly less than $4$.  
  The implementation \texttt{ISOP1.1}~\cite{Margot07} that calls the \texttt{CPLEX12.5.1} libraries was used for the B\&B with isomorphism pruning algorithm. Whenever isomorphism pruning was turned off, ILPs were solved by a B\&B algorithm implementation in {\tt CPLEX12.5.1} interactive. (The choice of the B\&B algorithm was made by optimizing the settings in \texttt{CPLEX} for $N=160,176$ on a small subset of  ILPs.) Jobs were run on our \texttt{Intel(R) Xeon(R)} $3.10$GHz and $2.00$GHz processors for the $N=160$ and $N=176$ cases. The $N=176$ cases were parallelized, but the $N=160$ cases were not.
 The MCS times in Table~\ref{tab:Schoen} are those reported by Schoen, Eendebak, and Nguyen~\cite{Schoen}.

An issue was determining the most efficient method for the computations in this paper.
  The method Full was the fastest for small $k \le t+3$. In these relatively quick classifications with small $k \le t+3$, a modest increase in real time was incurred (about double) if the Hybrid  method was used instead. While the classifications presented in Table~\ref{tab:enum} were still possible without using the method Full, it is included since the derivation of the  Hybrid  method requires an understanding of the full ILP formulation of the method Full. For the more time consuming problems with moderate and large $k > t+3$, we typically tested each method on the same small subset of $\oankm$ and then used the fastest to finish the full job. The Hybrid  method was the fastest with moderate $k=t+4$ and was the fastest by orders of magnitude when $k=t+4$ and  $N=160,176$. With large $k \ge t+5$, on the other hand, the ILP-extension method without isomorphism pruning  was the fastest and was the fastest by orders of magnitude when $N=24,28,36$.

As $N$ gets larger compared to $k$, the MCS algorithm suffers from the exponential growth of the search space and cannot compete with the best ILP-based method. However, when $N$ is small and the number of non-isomorphic $\oankm$ is large, the MCS algorithm is the fastest; see Table~\ref{tab:Schoen}. This happens for two reasons. (i) As the number of variables increases, solving LP relaxations can detect infeasibility higher in the search tree compared to the consistency checks of the MCS  algorithm.  
This makes the computational burden of solving LP relaxations worthwhile. (ii) The MCS algorithm does not find all non-isomorphic extensions of each non-isomorphic $\oankm$ in $M_{\ell}$ to $\oan$, but only those that are lexicographically minimum in columns. This eliminates many OAs that are otherwise produced by the ILP-based methods. Using the group $G_{s,k}$ in the Hybrid  method would get rid of this inefficiency.  However, using $G_{s,k}$ causes isomorphism classes of $\oan$ to  be lost because not all $\oan$ that are lexicographically minimum in rows under the action of $G_{s,k}$ are extensions of $\oankm$ possessing the same property under the action of $G_{s,k-1}$.

The times in Table~\ref{tab:Schoen} for ILP formulations are from  an earlier implementation of the B\&B with isomorphism pruning algorithm~\cite{Margot07}. This is built on
\texttt{Bcp/stable/1.1}~and uses \texttt{CLP1.7}; see \url{https://www.coin-or.org/wordpress/coin-xml/XML/Clp.xml} for solving LP relaxations.
We observed this implementation to be more than $6$ times slower than \texttt{CPLEX12.5.1} when its  isomorphism pruning  utility was turned off, so the ILP-based methods' times in Table~\ref{tab:Schoen} are conservatively high. In spite of this, the superiority of ILP-based methods over MCS at large $N$ is still the trend.

The MCS algorithm of Schoen, Eendebak, and Nguyen~\cite{Schoen} was compiled on our computer, and the resulting executable was used to directly address 3 fundamental questions. (i) Could $\text{OA}(160,k,2,4)$ classifications be completed by the  MCS algorithm together with OD-equivalence reductions? To answer this first question, $1,000$  randomly chosen $\text{OA}(160,8,2,4)$ were extended to $\text{OA}(160,9,2,4)$ using MCS and ILP-extension. ILP-extension was more than $3$ orders of magnitude faster, so ILP-extension was essential in this classification. (ii) Was isomorphism pruning essential in our ILP-based methods? To answer this second question, the isomorphism pruning option of  Hybrid  was turned off. This resulted in a 1 order of magnitude slow down for the extensions to $\text{OA}(160,8,2,4)$, and Step~\ref{item:nauty} of the Basic Extension method (\texttt{nauty} step) in extending to $\text{OA}(176,8,2,4)$ failed because the $64$GB of disk space mounted on our  \texttt{/tmp} directory was filled. (iii) Could the new cases be classified with only   \texttt{CPLEX12.5.1}  and the MCS algorithm, without isomorphism pruning? To answer this last question, the MCS algorithm was tested on cases where ILP-based methods required isomorphism pruning. The extension of all non-OD-equivalent $\text{OA}(160,7,2,4)$ to all non-OD-equivalent $\text{OA}(160,8,2,4)$ was $4.64$ times slower with the MCS algorithm. Furthermore, extending $\text{OA}(176,7,2,4)$ to $\text{OA}(176,8,2,4)$ was estimated to be about $11$ times slower based on $20$ randomly selected non-isomorphic input 
$\text{OA}(176,7,2,4)$s. Hence, these computational experiments suggest that   B\&B with isomorphism pruning provides a substantially larger magnitude of speed up   as $N$ increases. 

McKay~\cite{McKayfree,Radiz}  developed a technique for isomorph rejection for constructing combinatorial objects in which an inductive process of extending smaller objects to bigger ones as in Method~\ref{alg:basic}  exists, see Section 4.2.3 in~\cite{Kaski06}. McKay's technique removes some of the  isomorphisms only due to permutation of rows or columns. So, McKay's technique generates isomorphic (OD-equivalent) OAs.   We implemented McKay's technique within the method Hybrid by replacing  B\&B with isomorphism pruning with our choice of a basic 
B\&B implementation in {\tt CPLEX12.5.1} with McKay's technique for isomorph rejection. We removed the remaining isomorphic (OD-equivalent) OAs by using {\tt nauty} as in the Basic Extension method. 
 We tested the resulting method in extending all non-OD-equivalent OA$(N,7,2,4)$ to all non-OD-equivalent 
OA$(N,8,2,4)$ for $N=160,176$, i.e., only the cases in which isomorphism pruning 
was necessary. B\&B with isomorphism pruning was $11.44$ and $1.44$ times faster for the $N=160, 176$ cases.  The improvement  was much more dramatic when $N=160$ because there are many more exploitable symmetries by isomorphism pruning. (For the $N=160$ and $N=176$ cases the average sizes of the groups of all exploitable symmetries per input OA$(N,7,2,4)$ are $43.38$ and $3.11$.)  
We expect that for $N=192$ the improvement will be even more dramatic than in the $N=160$ case because when $N=192$ the average number of exploitable symmetries per input OA$(192,7,2,4)$ is $481$.
For enumerating all non-isomorphic $\text{OA}(192,k,2,4)$, our results suggest the use of OD-equivalence reductions with the method: Full when $k=5,6,7$, Hybrid when $k=8$, and ILP-extension when $k>8$.
\section{Consistency checks}\label{sec:validation}
 \setcounter{dummyt}{0} 
\setcounter{dummym}{0} 
\setcounter{dummya}{0} 
\setcounter{dummyc}{0}
\setcounter{dummyl}{0}
\setcounter{dummyd}{0}
\setcounter{dummye}{0}
Let $\mathcal{OD}_{N,k,t}$ be the set of all arrays that are Hadamard equivalent to $[\1,\C]$, where $\C$ is an $\oant$ with  symbols $\pm 1$. For $\Y_1,\Y_2 \in \mathcal{OD}_{N,k,t}$, let $\Y_1$ be considered to be equal to $\Y_2$ if $\Y_1$ can be obtained from $\Y_2$ by applying a signed permutation to rows. Hence, $|\mathcal{OD}_{N,k,t}|$ is the number of $[\1, \C]$ up to signed permutations of rows, where $\C$ is an $\text{OA}(N,k,2,t)$. 
The following lemma shows that $|\mathcal{OD}_{N,k,t}|$ is also the number of all $\oant$ up to permutations of rows.
\begin{lem}\label{signperm}
If $\C_1$ and $\C_2$ are $2$-symbol designs with $k$ columns and $N$ rows, then $\C_1=\C_2$ up to a permutation of rows if and only if $[\1 , \C_1]=[\1 , \C_2]$ up to a signed permutation of rows. Hence, $|\mathcal{OD}_{N,k,t}|$ is the number of all $\oant$ up to permutations of rows.
 \end{lem}
\begin{pf}
If $[\1 , \C_1]=[\1 , \C_2]$ up to a signed permutation of rows, then
$\vv_1\odot\pi_1([\1 , \C_1])=\vv_2\odot\pi_2([\1 , \C_2])$ for some $\pi_1,\pi_2 \in S_N$ and $\vv_1,\vv_2 \in \{\pm 1\}^N$. Also, $\vv_1=\mathbf \vv_2$ otherwise $\pi_1(\1)\neq \pi_2(\1)$. Hence, $\pi_1([\1,\C_1]) = \pi_2([\1,\C_2])\Rightarrow \pi_1(\C_1) = \pi_2(\C_2)$.
The converse is obvious. \qed
\end{pf}

Next, the method from Section 10.2 of~\cite{Kaski06} is applied to directly validate the correctness of the OD classifications summarized in Table~\ref{tab:enum}  by computing $|\mathcal{OD}_{N,k,t}|$ twice, i.e., once from the set of all non-OD-equivalent $\text{OA}(N,k,2,4)$ obtained by our methods and the second time based on that of $\text{OA}(N,k-1,2,4)$. 

The group of signed permutations $S_2\wr S_{k+1}$ of columns of matrices with $ (k+1)$ columns  acts on the elements of   $\mathcal{OD}_{N,k,t}$. Let 
$$\text{Aut}([\1 , \C])=\{\pi \in S_2\wr S_{k+1}\ |\  \sigma [\1 , \C]\pi =[\1 , \C] \text{ for some } \sigma \in S_2\wr S_N\},$$
where $S_2\wr S_N$ is the group of signed permutations of rows of matrices with $N$ rows.
Then by the orbit-stabilizer theorem, when $\C$ is an $\oant$, the total number of elements in $\mathcal{OD}_{N,k,t}$ that are Hadamard equivalent to the element  $[\1,\C] \in \mathcal{OD}_{N,k,t}$ is $$\frac{(k+1)!2^{k+1}}{|\text{Aut}([\1 , \C])|}.$$
 Hence,
\begin{equation}\label{eqn:count1}
|\mathcal{OD}_{N,k,t}|=\sum_{[\1, \C] \+ \in \+ \reallywidehat{\CC}_{N,k,t}}\frac{(k+1)!2^{k+1}}{|\text{Aut}([\1, \C])|}=
\sum_{\C \+ \in \+ \mathlarger{\CC}_{N,k,t}}\frac{(k+1)!2^{k+1}}{|\text{Aut}([\1, \C])|},
\end{equation}
where $\reallywidehat{\CC}_{N,k,t}$ is a set of all 
non-Hadamard equivalent elements in $\mathcal{OD}_{N,k,t}$ and $\CC_{N,k,t}$ is a set of all non-OD-equivalent $\oant$. Clearly, $|\CC_{N,k,t}|=|\reallywidehat{\CC}_{N,k,t}|$.
(For example, from Table~\ref{tab:enum} we have $|\CC_{160,8,4}|=|\reallywidehat{\CC}_{160,8,4}|=11,712$, and the sum in~(\ref{eqn:count1}) for $|\mathcal{OD}_{160,8,4}|$ has $11,712$ terms.) 

Also, for fixed $\C'\in \CC_{N,k-1,t}$ let $N([\1,\C'])$ be the number of distinct $[\1,\C] \in \mathcal{OD}_{N,k,t}$ that can be obtained by adding a column to $[\1,\C'] \in \mathcal{OD}_{N,k-1,t}$. By the orbit-stabilizer theorem, the number of elements in $\mathcal{OD}_{N,k,t}$ whose first $k$ columns subarray is Hadamard equivalent to $[\1, \C']$ is $k!2^{k}N([\1,\C'])/|\text{Aut}([\1 , \C'])|$, so
\begin{equation}\label{eqn:count2}
|\mathcal{OD}_{N,k,t}|=\sum_{[\1, \C']\+ \in\+ \reallywidehat{\CC}_{N,k-1,t}}\frac{k!2^{k}N([\1, \C'])}{|\text{Aut}([\1, \C'])|}=\sum_{\C'\+ \in \+ \mathlarger{\CC}_{N,k-1,t}}\frac{k!2^{k}N([\1,\C'])}{|\text{Aut}([\1 , \C'])|}.
\end{equation}
For each $\mathbf{C}'$, the value of $N([\1, \mathbf{C}'])$ in~(\ref{eqn:count2}) is calculated by taking the OD-equivalence class representative $\mathbf{C}'$  as the input OA$(N,k-1,2,t)$ in  ILP-extension or  Hybrid and counting all possible solutions to ILP~(\ref{ilp:oa}) or ILP~(\ref{ilp:oahybrid}). This is done by finding all possible solutions using our choice of a basic B\&B implementation in {\tt CPLEX12.5.1}.
The counts from equations~(\ref{eqn:count1}) and~(\ref{eqn:count2})   must match as we are calculating  $|\mathcal{OD}_{N,k,t}|$ in two different ways, once using the set of all non-OD-equivalent OA$(N,k,2,t)$ in equation~(\ref{eqn:count1}) and a second time using the set of all  non-OD-equivalent OA$(N,k-1,2,t)$ in equation~(\ref{eqn:count2}). However, a mismatch is possible as we are using two different equations to calculate $|\mathcal{OD}_{N,k,t}|$. 
If we do not get the same number  $|\mathcal{OD}_{N,k,t}|$ using the two different equations, then there may be missing OD-equivalence classes of OA$(N,k-1,2,t)$ or OA$(N,k,2,t)$, or at least one of the automorphism group sizes might have been incorrectly calculated.
 
Recall by Lemma~\ref{signperm} that $|\mathcal{OD}_{N,k,t}|$ is also equal to the number of $\text{OA}(N,k,2,t)$ up to  permutations of rows. Then the numbers of $\text{OA}(160,k,2,4)$ up to permutations of rows obtained by both counts~(\ref{eqn:count1}) and~(\ref{eqn:count2}) are $11$, $5482$, $61084192$, $855133483008$, $1435727462400$  for $k=5,6,7,8,9$. Similarly, the numbers of  $\text{OA}(176,k,2,4)$ up to permutations of rows obtained by both counts~(\ref{eqn:count1}) and~(\ref{eqn:count2}) are $12$, $7680$, $400934400$, $11928390727680$, $4830658560$  for $k=5,6,7,8,9$. The automorphism group sizes  $|\text{Aut}([\1 , \C])|$ and $|\text{Aut}([\1 , \C'])|$ in sums~(\ref{eqn:count1}) and~(\ref{eqn:count2}) were computed by using {\tt nauty}~\cite{McKayP, McKay13}, and the sums~(\ref{eqn:count1}) and~(\ref{eqn:count2}) were calculated using 
{\tt GAP}~\cite{GAP} with infinite-precision integers to avoid errors due to limited precision. Hence the OD-equivalence class  sizes in Table~\ref{tab:enum} are validated.
Since a set of all non-isomorphic OAs was obtained from a set of all non-Hadamard equivalent elements of $\mathcal{OD}_{N,k,t}$ by the method in Section~\ref{sec:fold}, this consistency check also validates the OA classifications.

\section{Discussion and challenge problems}\label{sec:sum}
New ILP-based methods with the isomorphism pruning  option classified $\text{OA}(N,k,2,4)$ with $N=160,176$. A comparison of these methods to  the MCS algorithm~\cite{Schoen}, the best known algorithm for this purpose, was also made. The MCS algorithm does not solve LP relaxations in its search tree. Our paper clearly demonstrates the value of solving LP relaxations in extending OA$(N,k-1,2,t)$ to OA$(N,k,2,t)$ for $N \geq 160$ and $k\geq 9$.
It also shows the value of combining  dynamic isomorph rejection with solving LP relaxations in a B\&B tree for $N \geq 160$ and  $k=8$. This is made possible by the B\&B with isomorphism pruning algorithm~\cite{Margot07}.

We have no method of reducing the number of ILPs that need to be solved analogous to OD-equivalence reductions when $s\geq3$. Consequently, we have not been able to classify all isomorphism classes of $\oan$ for  new cases  when $s\geq3$. However, as documented in Table~\ref{tab:Schoen}, our methods classified all non-isomorphic $\text{OA}(81,k,3,3)$ much faster than the MCS algorithm~\cite{Schoen}, but this speed improvement is not sufficient to classify all isomorphism classes of $\text{OA}(108,k,3,3)$.

Our methods become less effective with increased $s$ and $k$. This is due to exponentially increasing sizes of the exploitable symmetry groups and the number of variables. In particular, methods in this paper are not competitive with those in~\cite{Egan,Kokkala15} for the OAs they classified.

Using the method Full we found that there are $1,529$ OD-equivalence classes and $8,386$ isomorphism classes of  OA$(192,7,2,4)$ in $174$ minutes real time
on a $3.10$ GHz processor. 
However, we could not find a set of all non-isomorphic $\text{OA}(192,k,2,4)$ for $k\geq8$ because the number of non-OD-equivalent $\text{OA}(192,8,2,4)$ is too large. We propose determining $f(192,2,4)$, $f(108,3,3)$, and $f(100,10,2)$ as challenge problems. Among these, the solution to the problem involving  OA$(100, k, 10, 2)$ is the most sought after in  combinatorics. We do not even know whether  an OA$(100, 5, 10, 2)$ exists.
 
 In~\cite{AMM}, it is shown that embedding a B\&B algorithm to solve the ILP formulation in~\cite{BulutogluM} within a constraint programming
 search tree is faster than using the B\&B alone on the same formulation in determining the existence of  an OA$(q^2,4,q,2)$  for $3 \leq q \leq 12$.
 Hence, for classifying $\oan$, a constraint programming with isomorphism pruning algorithm that  solves LP relaxations, generates cutting planes, and resolves  to prune  by infeasibility only at certain depths of the search tree might be beneficial.

Developing two different methods based on two different ILP formulations of the orthogonal array problem was essential
in bringing the classification of  $\text{OA}(N,k,2,4)$ for $N=160,176$ within computational reach. 
  Determining formulations better than the full ILP 
 formulation~(\ref{ilp:BF}) that also do not require enumerating all non-OD-equivalent $\oan$  is another direction to look at for solving the aforementioned challenge problems.
  
 A collection of $b=\lambda {\upsilon \choose t}/ {k \choose t}$
size $k$ sets (blocks) of $\upsilon$ treatments (symbols)  such that each of the ${\upsilon \choose t}$ possible $t$ treatment combinations appears in  $\lambda$ blocks is called a 
{\em $t$-$(\upsilon,k,\lambda)$ design}. Two $t$-$(\upsilon,k,\lambda)$ designs are {\em isomorphic} if one can be obtained from the other up to a permutation of blocks by permuting the treatments.
An $\oan$ and a $t$-$(\upsilon,k,\lambda)$ design   are analogous combinatorial objects as they  are both $t$-designs in the Hamming and Johnson cometric schemes respectively, see  Chapter 30 and Theorem 30.5 in~\cite{Wilson} for the definitions and the result. Hence, there are analogous ILP formulations for classifying $t$-$(\upsilon,k,\lambda)$ designs up to isomorphism to those for classifying $\oan$. In particular, ILP~(2.1) in~\cite{Margot03a} after replacing ``$\geq$" with ``$=$"  is a formulation for classifying $t$-$(\upsilon,k,\lambda)$ designs   that is analogous to the Bulutoglu and Margot~\cite{BulutogluM} ILP formulation for classifying OAs. This formulation leads to a method analogous to the method Full.
The blocks, treatments, the parameter $t$, and the parameter $\lambda$ of a $t$-$(\upsilon,k,\lambda)$ design are analogous to the rows, columns, the strength, and the index ($N/s^t$) of an $\oan$. So, an ILP-based method analogous to ILP-extension or Hybrid based on finding all possible ways of adding a new treatment to each isomorphism class representative of all partial $t$-$(\upsilon,k,\lambda)$ designs with $i-1$ treatments to get a set of all non-isomorphic partial $t$-$(\upsilon,k,\lambda)$ designs with $i$ treatments for $i=t+1,\ldots,\upsilon$ can be formulated. This is the ILP-based version of the column-by-column construction (ALG1) in~\cite{Radiz}. Alternatively, the ILP formulation within the method analogous to 
Full can be used to develop another ILP formulation that finds a set of all $(t'+1)$-$(\upsilon'+1,k'+1,\lambda)$ designs by extending a $t'$-$(\upsilon',k',\lambda)$ design with a new treatment. Then, this formulation can be used to develop a method that finds a set of all non-isomorphic $(t'+1)$-$(\upsilon'+1,k'+1,\lambda)$ designs by extending a set of all non-isomorphic $t'$-$(\upsilon',k',\lambda)$ designs with a new treatment.
Such a method is the ILP-based version of the row-by-row construction (ALG2) in~\cite{Radiz}.
For classifying  $t$-$(\upsilon,k,\lambda)$ 
designs up to isomorphism, we expect the observed trend from Table~\ref{tab:Schoen} to hold between the aforementioned 
ILP-based methods analogous to those in this paper and backtrack algorithms with isomorph rejection that do not solve LP relaxations within their search trees.  Hence,  B\&B  with isomorphism pruning~\cite{Margot07} or a basic B\&B  coupled with McKay's technique from~\cite{Radiz}    for isomorph rejection may bring some unsolved cases of classifying  $t$-$(\upsilon,k,\lambda)$ 
designs for  $\lambda \geq 2$ and larger values of $b=\lambda {\upsilon \choose t}/ {k \choose t}$ within computational reach. The smallest open problems of classifying all non-isomorphic $4$-$(\upsilon,k,\lambda)$ are $4$-$(11,5,3)$, $4$-$(12,6,8)$, $4$-$(13,6,6)$, $4$-$(13,6,12)$, $4$-$(13,6,18)$, $4$-$(15,5,2)$, and $4$-$(18,6,2)$ designs~\cite{Kaski06,Radiz}.
\section*{Acknowledgements}
The authors  thank an anonymous referee for improving Section~\ref{sec:validation}.  The authors also thank Professor~Fran\c cois~Margot for providing his implementation of the B\&B with isomorphism pruning algorithm and answering questions regarding his software; Professor~Brendan McKay for answering questions regarding his program \texttt{nauty}; and Mr. David Doak for general computer support.

 This research was supported by the AFOSR  grants F1ATA06334J001 and F1ATA03039J \ 001. The views expressed in this article are those of the authors and do not reflect the official policy or position of the United States Air Force, Department of Defense, or the U.S.~Government.   This paper has been published at the Australasian Journal of Combinatorics: 
  \url{http://ajc.maths.uq.edu.au}.
\begin{appendix}{}
\section{Tables}\label{sec:tables}
\scriptsize
\begin{table}[h]
\caption{The number of OA isomorphism classes, the number of OD-equivalence classes,
real time (minutes), method.
} \label{tab:enum} \vspace{.1in} \centering
\begin{tabular}{|lrrrc|}\hline
$N,k,s,t$ & OA & OD & Time & Method  \\ \hline
 $160,5,2,4$ &        6 &  &   0.033 &      Full  \\
 $160,6,2,4$ &       29 &  &   0.5 &      Full  \\
 $160,7,2,4$ &      450 &106&  4.5 &      Full  \\
 $160,8,2,4$ &    99,618 &11,712& 1,184 &  Hybrid  \\
 $160,9,2,4$ &    15,083 &1,608&4,384 &  ILP-Ext. \\
$160,10,2,4$ &        0 &0 &   4.1 &   ILP-Ext.   \\ \hline
 $176,5,2,4$ &        6 &  &  0.01 &      Full  \\
 $176,6,2,4$ &       14 &  &  0.09 &      Full  \\
 $176,7,2,4$ &      945 &179&  15.3 &  Full   \\
 $176,8,2,4$ &  1,157,443 &129,138& 32,470 & Hybrid \\
  $176,9,2,4$ &  26 &  4& 378,241 & ILP-Ext. \\
  $176,10,2,4$ &  0 &  0& 0.01 & ILP-Ext. \\ \hline
\end{tabular}
\end{table}

\begin{table}[h]
\caption{Real times (minutes) of OA classifications using our ILP methods without OD-equivalence reductions and the MCS algorithm of Schoen, Eendebak, and Nguyen~\cite{Schoen}.} \label{tab:Schoen} \vspace{.1in} \centering
\begin{tabular}{|l|rr|}\hline
$\oan$ & ILP & MCS \\ \hline
$\text{OA}(24,k,2,2)$  & 1,396,589 & 66,369 \\
$\text{OA}(32,k,2,3)$  &       1.2 &    4.9 \\
$\text{OA}(40,k,2,3)$  &        12 &    4.8 \\
$\text{OA}(48,k,2,3)$  &   605,765 & 74,232 \\
$\text{OA}(81,k,3,3)$  &     2,610 &      31,660 \\ \hline
$\text{OA}(96,k,2,4)$  &       0.4 &   0.02 \\
$\text{OA}(112,k,2,4)$ &       0.1 &   0.02 \\
$\text{OA}(128,k,2,4)$ &       914 &  4,372 \\
$\text{OA}(144,k,2,4)$ &        37 &    229 \\
$\text{OA}(160,k,2,4)$ &   711,298 &      ? \\ \hline
\end{tabular}
\end{table}

\end{appendix}
\vspace{3cm}


\end{document}